\documentclass[12pt]{amsart}
\usepackage{amssymb,latexsym}
\usepackage{enumerate}
\usepackage[mathscr]{eucal}
\usepackage[isolatin]{inputenc}
\usepackage{upref}
\usepackage{esint}
\usepackage{epsfig}
\usepackage{verbatim}
\usepackage{multicol}
\usepackage{multirow}
\usepackage{graphicx, color, psfrag, epstopdf}
\usepackage{pstricks, pst-node}

\usepackage[paper=a4paper,left=30mm,right=20mm,top=25mm,bottom=30mm]{geometry}

\makeatletter
\@namedef{subjclassname@2010}{
 \textup{2010} Mathematics Subject Classification}
\makeatother

\newtheorem{thm}{Theorem}[section]

\newtheorem{lem}[thm]{Lemma}

\theoremstyle{definition}
\newtheorem{defin}[thm]{Definition}
\newtheorem{rem}[thm]{Remark}

\numberwithin{equation}{section}

 \newcommand{\N}{\mathbb{N}}
 \newcommand{\R}{\mathbb{R}}

\newcommand{\dokendDef}{\ensuremath{\hfill\square}}

\begin{document}
%%%%%%%%%%%%%%%%%%%%%%%%%%%%%%%%%%%%%%%%%%%%%%%%%%%%%%%%%%%%%%%%%%%%%%%
% Mit LaTeX=>PS=>PDF uebersetzen, damit Figure 1 generiert werden kann!
%%%%%%%%%%%%%%%%%%%%%%%%%%%%%%%%%%%%%%%%%%%%%%%%%%%%%%%%%%%%%%%%%%%%%%%
\baselineskip=17pt

\title[Radially symmetric solutions]
{Radially symmetric solutions of the ultra-relativistic Euler equations}

\author{Matthias Kunik}
\address{Otto-von-Guericke-Universit\"{a}t Magdeburg\\
Institut f\"ur Analysis und Numerik\\
Geb\"{a}ude 02 \\
Universit\"{a}tsplatz 2 \\
D-39106 Magdeburg \\
Germany}
\email{matthias.kunik@ovgu.de}

\author{Hailiang Liu}
\address{Iowa State University\\ Department of Mathematics\\ Ames, IA 50010-2064, USA}
\email{hliu@iastate.edu}

\author{Gerald Warnecke}
\address{Otto-von-Guericke-Universit\"{a}t Magdeburg\\
Institut f\"ur Analysis und Numerik\\
Geb\"{a}ude 02 \\
Universit\"{a}tsplatz 2 \\
D-39106 Magdeburg \\
Germany}
\email{gerald.warnecke@ovgu.de}

\date{\today}

%\begin{dedication}
%Dedicated to Professor Hsiao Ling on the occasion of her 80th birthday
%\end{dedication}
\dedicatory{Dedicated to Professor Hsiao Ling on the occasion of her 80th birthday}

\begin{abstract}
The ultra-relativistic Euler equations for an ideal gas are described in terms of the pressure $p$,
the spatial part $\underline{u} \in \R^3$ of the dimensionless four-velocity
and the particle density $n$. Radially symmetric solutions of these equations are studied.
Analytical solutions are presented for the linearized system. For the original nonlinear equations we design and analyze a numerical scheme for simulating radially symmetric solutions in three space dimensions. The good performance of the scheme is demonstrated by  numerical examples. In particular, it was observed that the method has the capability to capture accurately the pressure singularity formation caused by shock wave reflections at the origin.
\end{abstract}



\maketitle

\section{Introduction}
Relativistic flow problems are vital in many astrophysical phenomena. 
An effective way to improve our knowledge of the actual mechanisms 
is due to relativistic hydrodynamics simulations.
Especially, solutions describing radially symmetric gas flow are important in applications as well as in theory.
They are particularly well suited for numerical simulations of certain multi-dimensional problems. In this paper 
we focus on radially symmetric solutions. We consider a special relativistic system which is much simpler than flows in general relativistic theory. Interestingly, even compared to the classical Euler equations of non-relativistic gas dynamics the equations we consider exhibit a simpler mathematical structure.

We are concerned with the ultra-relativistic equations 
for a perfect fluid in Minkowski space-time, namely
\begin{equation}\label{div}
\sum_{\beta=0}^{3}\frac{\partial T_{\alpha \beta}}{\partial x_{\beta}} =
0,~~~ \frac{\partial N_{\alpha}}{\partial x_{\alpha}} = 0,
\end{equation}
where
\begin{equation*}
T_{\alpha\beta} =-p g_{\alpha\beta}+4pu_\alpha u_\beta
\end{equation*}
denotes the energy-momentum tensor for the ideal ultra-relativistic gas.
Here $p$ represents the pressure,
$\underline{u} \in \R^3$ is the spatial part of the four-velocity
$(u_0,u_1,u_2,u_3)=(\sqrt{1+|\underline{u}|^2},\underline{u})$\,.
The flat Minkowski metric is given as
\begin{equation*}
g_{\alpha \beta}= \left\{\begin{array}{ccc}
 +1,   &\alpha =\beta= 0\,,\\
-1,   &\alpha =\beta = 1,2,3\,,\\
0,    & \alpha \neq \beta\,,
\end{array}
\right.
\end{equation*}
and the particle-density four-vector is denoted by 
\begin{equation}\label{Nalpha}
N_{\alpha}= n u_{\alpha}.
\end{equation}
Here $n$ is the proper particle density. 
We note that the quantities $u_{\alpha}$, $T_{\alpha \beta}$, $g_{\alpha \beta}$,
$N_{\alpha}$\, and even $x_{\alpha}$\, are usually written down as Lorentz-invariant
tensors with \textit{upper indices} instead of  lower indices in order to make use of
Einstein's summation convention. But in the following calculations these upper indices could be mixed up with powers.
Since we will not make use of the lowering and raising of Lorentz-tensor indices,
our change of the notation will not lead to confusions.
For the physical background we refer to Weinberg \cite[Part one, pp 47-52]{Weinberg},
further details can be found in Kunik \cite[Chapter 3.9]{Kunikthesis},
and for the corresponding classical Euler equations see Courant and Friedrichs \cite{CF}. For a general introduction to the mathematical theory of hyperbolic conservation laws see Bressan \cite{Bressan} and Dafermos \cite{CD}.  
A nice overview of radially symmetric solutions to conservation laws
is given in Jenssen's survey paper \cite{HJ}.

The unknown quantities $p$, $\underline{u}$
and $n$ satisfying \eqref{div} depend in general on
time $t=x_0 \geq 0$ and position $\underline{x}=(x_1,x_2,x_3) \in \R^3$\,.
It is well known that even for smooth initial data, where the fields 
are prescribed at $t=0$, the solution may develop shock discontinuities.
This requires a weak form of the conservation laws in \eqref{div}.
Since the conservation law for the particle-density four-vector \eqref{Nalpha} decouples from the conservation laws of 
energy and momentum, we will restrict ourselves to the resulting
closed subsystem for the variables $p$ and $\underline{u}$ 
satisfying the first set of equations in \eqref{div}. 
Putting $\alpha=0$ this
gives the conservation of energy
\begin{equation}\label{energy_general}
\frac{\partial}{\partial t}\left(3p+4p|\underline{u}|^2 \right)+
\sum\limits_{k=1}^{3}\frac{\partial}{\partial x_k}
\left(4pu_k\sqrt{1+|\underline{u}|^2} \right)=0\,,
\end{equation}
whereas for $\alpha=j=1,2,3$ we obtain the conservation of momentum
\begin{equation}\label{momentum_general}
\frac{\partial}{\partial t}\left(4pu_j \sqrt{1+|\underline{u}|^2} \right)+
\sum\limits_{k=1}^{3}\frac{\partial}{\partial x_k}
\left(p\delta_{jk}+4pu_j u_k\right)=0\,,\quad j=1,2,3\,.
\end{equation}
Like the classical Euler equations, these relativistic Euler equations constitute a hyperbolic system of conservation laws and have their origin in the kinetic theory of gases. This can be used for the construction of numerical schemes which preserve positive pressure and satisfy a discrete version of the entropy inequality, 
see \cite{Kunikthesis,KQW1,KQW3}. Some other analytical and numerical methods 
for the ultra relativistic Euler equations are studied in \cite{Abdelrahman1,
Abdelrahman2, Abdelrahman3, HK, KQW2,  QYM11, QY12, QY13, HT12}. Recently numerical results using central upwind scheme are reported  in \cite{GYQ18} for one and two dimensional special ultra-relativistic Euler equations.
In \cite{GL} Lai presents a detailed analysis of self-similar solutions of 
radially symmetric relativistic Euler equations in three and two space dimensions.
These are special solutions depending only on $r/t$ with radius $r$ and time $t$
which satisfy systems of ordinary differential equations. Especially his study of the ultra-
relativistic Euler equations enables us to compare his solutions with two of our numerical
results.

In this paper, we study radially symmetric solutions and construct a corresponding scheme to solve the ultra-relativistic Euler equations \eqref{energy_general}, \eqref{momentum_general} in three space dimensions. One of the main advantages of the radially symmetric problem 
is that it can be used to efficiently simulate special wave patterns for fully three dimensional problems such as the 
detonation problem, see \cite{To99}; also Example 5.5 in \cite{SL11}
for the classical Euler equations. We show this with Example 4 
in Section \ref{numeric}
for the ultra-relativistic Euler equations. This allows the prediction of pressure singularity formation caused by shock wave reflections at the origin. The shock wave reflection 
with a pressure singularity at the boundary is also motivated by the analysis of the linearized model in Section \ref{linear}.
We hope that our specific solutions may become benchmarks for testing fully 3D simulations.

In the next section we will define radially symmetric solutions of \eqref{energy_general}, \eqref{momentum_general}
in a weak integral form. We will also present the initial- and boundary value problem for nonlinear radially symmetric solutions. 
In Section \ref{linear} we will especially solve the corresponding linearized model,
which reveals singularities in the pressure field due to shock reflections at the boundary. 
Formulation and analysis of a stable scheme
based on the balance laws \eqref{weak_contour} 
are presented in Section \ref{scheme1}. In this connection emphasis is on a proper treatment of the radius 
and the boundary conditions in the balance laws. 
The numerical examples are given in Section \ref{numeric}.

\section{Radially symmetric solutions}\label{radial_solutions}
Assume for a moment a smooth solution $p, \underline{u}$ of the ultra-relativistic Euler equations 
\eqref{energy_general}, \eqref{momentum_general}. We put $r = |\underline{x}|$ for $r\geq 0$
and look for radially symmetric solutions
\begin{equation}\label{rad}
p=p(t,r)>0\,, \quad \underline{u}(t, \underline{x})=\frac{u(t,r)}{r}\,\underline{x}\,.
\end{equation}
Here the quantity $\underline{u}(t, \underline{x}) \in \R^3$ is completely determined
by a new \textit{real valued quantity} $u(t,r)$ depending on $t > 0$, $r > 0$.
For continuity we have the boundary condition
\begin{equation}\label{zeroboundary}
\lim \limits_{r \downarrow 0} u(t,r)=0\,, \quad t > 0\,.
\end{equation}
Note that $\underline{n} = \frac{1}{r}\underline{x}$ is the outer normal vector field 
of the sphere $\partial \mathcal{B}_R$ bounding the ball 
$\mathcal{B}_R=\{\underline{x} \in \R^3\,:\,|\underline{x}| \leq R \}$ of radius $R>0$,
and that $|\underline{u}|^2=u^2$ as well as $u=\underline{u}\cdot \underline{n}$\,.
Therefore, it is natural to apply the Gaussian divergence theorem for
the integration of the second term with respect to $\mathcal{B}_R$ of the conservation law 
\eqref{energy_general} in order to make use of the radial symmetry of the fields.
We obtain with \eqref{rad} for any fixed $R>0$
\begin{equation*}
4\pi \frac{\partial}{\partial t}\int \limits_0^R\left(3p(t,r)+4p(t,r)u^2(t,r)\right)\,r^2 dr
+\int \limits_{\partial \mathcal{B}_R}4pu\sqrt{1+u^2}\,dS=0\,.
\end{equation*}
The integrand in the surface integral is constant. Hence we have
\begin{equation}\label{energy_int}
\begin{split}
& \frac{\partial}{\partial t}\int \limits_0^R\left(3p(t,r)+4p(t,r)u^2(t,r)\right)\,r^2 dr\\
& +R^2 4p(t,R)u(t,R)\sqrt{1+u^2(t,R)}=0\,.\\
\end{split}
\end{equation}
This idea does not work for the momentum equation \eqref{momentum_general},
because \eqref{rad} would give values zero after integration with respect to $\mathcal{B}_R$.
Here we integrate \eqref{momentum_general} for $j=3$ over the upper half-ball
$$\mathcal{B}^+_R=\{\underline{x}=(x_1,x_2,x_3) \in \R^3\,:\,x_3 \geq 0\}\,,$$
use the Gaussian divergence theorem and spherical coordinates 
\begin{equation*}
x_1=r\cos{\varphi}\sin{\vartheta}\,, x_2=r\sin{\varphi}\sin{\vartheta}\,, x_3=r\cos{\vartheta}
\end{equation*}
with $r>0$, $0 < \varphi <2\pi$ and $0 < \vartheta < \pi/2$, and  obtain from \eqref{rad} 
\begin{equation}\label{momentum_int}
\begin{split}
 & \frac{\partial}{\partial t}\int \limits_0^R 4p(t,r)u(t,r)\sqrt{1+u^2(t,r)} \, r^2 dr\\
& +R^2\left(4p(t,R)u^2(t,R)+p(t,R)\right) =2\int \limits_0^R p(t,r)\,rdr.\\
\end{split}
\end{equation}

Now we differentiate the equations \eqref{energy_int}, \eqref{momentum_int} 
with respect to $R>0$. Afterwards we replace $R$ by the better suited variable $x>0$\,.\\
We put $p=p(t,x)$, $u=u(t,x)$ for abbreviation and have the $2$ by $2$ system
\begin{equation}\label{ultra_rad1}
\left\{
\begin{array}{ll}
 & \begin{displaystyle}\frac{\partial}{\partial t}\left(x^2 p (3+4u^2) \right)
+ \frac{\partial}{\partial x}\left(4 x^2 p u\sqrt{1+u^2} \right)\end{displaystyle}=0\,,\\[1.2ex]
& \begin{displaystyle}\frac{\partial}{\partial t}\left(4 x^2 p u\sqrt{1+u^2}  \right)
+ \frac{\partial}{\partial x}\left(x^2 p (1+4u^2) \right)\end{displaystyle}=2xp\,.\\
\end{array}
\right.
\end{equation}
The validity of this system may also be checked by differentiation 
from \eqref{rad}, \eqref{energy_general} and \eqref{momentum_general}.
The solutions of \eqref{ultra_rad1} are restricted to the state space
$\mathcal{S}_{eul}=\{(p,u) \in \R^2\,:\,p>0\}$.

For the formulation of weak entropy solutions we will introduce a transformation in state space. With 
$\begin{displaystyle}\tilde{\mathcal{S}}_{eul}=\{(a,b) \in \R^2\,:\,|b| < a\}\end{displaystyle}$
there is a one-to-one transformation $\Theta: \mathcal{S}_{eul} \mapsto \tilde{\mathcal{S}}_{eul}$
given by
\begin{equation}\label{statetrans}
\Theta(p,u)=
\begin{pmatrix}
p(3+4u^2)\\
4pu\sqrt{1+u^2}\\
\end{pmatrix}=\begin{pmatrix}
a\\
b\\
\end{pmatrix}\,.
\end{equation}
The inverse transformation is given by
\begin{equation*}
p=\frac{1}{3}\left( \sqrt{4a^2-3b^2}-a \right)\,,\quad
u=\frac{b}{\sqrt{4p(p+a)}}\,.
\end{equation*}
Using the transformation \eqref{statetrans} in state space we can also rewrite \eqref{ultra_rad1}
in an equivalent form. We put 
\begin{equation}\label{cdefinition}
c = c(a,b)=\frac{5}{3}a-\frac{2}{3}\sqrt{4a^2-3b^2}\,,
\end{equation}
and obtain from \eqref{ultra_rad1}, \eqref{statetrans}
\begin{equation}\label{ultra_rad3}
\left\{
\begin{array}{ll}
 & \begin{displaystyle}\frac{\partial}{\partial t}\left(x^2 a \right)
+ \frac{\partial}{\partial x}\left(x^2 b \right)\end{displaystyle}=0\,,\\[1.2ex]
& \begin{displaystyle}\frac{\partial}{\partial t}\left(x^2 b \right)
+ \frac{\partial}{\partial x}
\left(x^2 c \right)\end{displaystyle}
=x(a-c)\,.\\
\end{array}
\right.
\end{equation}
We look for weak solutions of \eqref{ultra_rad3} 
in the quarterplane $$Q=\{(t,x)\,:\,t>0\,,x>0\}\,.$$
For $x>0$ we prescribe the \textit{two initial functions}
\begin{equation}\label{initial_final}
\lim \limits_{t \,\downarrow \, 0}a(t,x)=a_0(x)\,, ~
\lim \limits_{t \,\downarrow \, 0}b(t,x)=b_0(x)\,, ~ x>0
\end{equation}
with $|b_0(x)|<a_0(x)$ for $x>0$.
Recall our preliminary assumption that we have a smooth and radially symmetric
solution of the three dimensional ultra-relativistic Euler equations \eqref{energy_general}
and \eqref{momentum_general}, which implies the boundary condition \eqref{zeroboundary}
as well as a locally bounded energy and momentum density.
Hence we have the following \textit{two boundary conditions} for $t>0$:
\begin{equation}\label{boundary_final}
\lim \limits_{x\,\downarrow\,0} (x^2 a(t,x))=0\,, \quad
\lim \limits_{x\,\downarrow\,0}b(t,x)=0\,.
\end{equation} 
We recall \eqref{cdefinition}, multiply in both equations of \eqref{ultra_rad3} with 
any $C^1$ test function $\phi:\R^2 \mapsto \R$ with compact support in $\R^2$
and obtain from \eqref{boundary_final} after partial integration 
\begin{equation}\label{ultra_rad2}
\left\{
\begin{array}{ll}
&\begin{displaystyle}\iint \limits_{Q} \left(a \, \frac{\partial \phi}{\partial t}
 + b\, \frac{\partial \phi}{\partial x} \right)\,x^2 dx dt
 +\int \limits_{0}^{\infty}\,a_0 \, \phi_0 \, x^2 dx
\end{displaystyle}=0\,,\\[1.4ex]
&\quad \begin{displaystyle}\iint \limits_{Q} \left(
 b\, \frac{\partial \phi}{\partial t} 
 +c(a,b)\,\frac{\partial \phi}{\partial x} \right)\, x^2 dx dt\end{displaystyle} ~ \\[1.2ex]
& \begin{displaystyle}+\iint \limits_{Q}(a-c(a,b))\,\phi \, xdx dt +
\int \limits_{0}^{\infty} b_0\, \phi_0 \, x^2 dx\end{displaystyle} = 0\,.\\
\end{array}
\right.
\end{equation}
We use $\phi=\phi(t,x)$, $\frac{\partial \phi}{\partial t}=\frac{\partial \phi}{\partial t}(t,x)$,
$\frac{\partial \phi}{\partial x}=\frac{\partial \phi}{\partial x}(t,x)$, $a=a(t,x)$, $b=b(t,x)$, 
$\phi_0=\phi(0,x)$, $a_0=a_0(x)$, $b_0=b_0(x)$ 
as abbreviations in \eqref{ultra_rad2}.
Now we drop the assumption that we have a smooth solution of the ultra-relativistic Euler equations
and will no longer assume that $a$ and $b$ are locally bounded.
\begin{defin}\label{def:weakrad} {\bf Weak radially symmetric solutions}\\
We say that $a,b$ is a \textit{weak solution}
of \eqref{ultra_rad3} with initial data $a_0,a_0$ if and only if 
the following conditions are satisfied:
\begin{itemize}
\item $a,b: Q \mapsto \R$ are measurable with $|b|<a$\,.
\item $x a(t,x)$ is integrable in $(0,t_0)\times(0,x_0) \subset Q$
for all $t_0,x_0>0$\,.
\item $a_0,b_0 : \R_{>0} \mapsto \R$ are measurable with $|b_0|<a_0$. We require for all $x_0>0$
that $x a_0(x)$ is integrable for $0<x<x_0$.
\item The boundary conditions \eqref{boundary_final} are satisfied
for almost all $t>0$.
\item Equations \eqref{ultra_rad2} are satisfied for all
 $C^1$ test function $\phi:\R^2 \mapsto \R$ with compact support in $\R^2$.
\end{itemize}
If $\psi : Q \mapsto \R_{\geq 0}$ is a new \textit{nonnegative} $C^1$ test function restricted to the 
quarter plane $Q$ with compact support in $Q$, then we will consider a weak solution $a$, $b$
which further satisfies
the \textit{weak entropy inequality}, see Kunik \cite[Chapter 4.4]{Kunikthesis},
\begin{equation}\label{weak_entropy}
\iint \limits_{Q} \left(p^{3/4}\sqrt{1+u^2} \, \frac{\partial \psi}{\partial t}
 + p^{3/4} u \, \frac{\partial \psi}{\partial x} \right)\,x^2 dx dt  \leq 0\,.
\end{equation}
Here we make use of $\Theta(p,u)=\begin{pmatrix}a\\b\\ \end{pmatrix}$
in \eqref{statetrans}.
In this case we call $a,b$ a \textit{weak entropy solution} of the system \eqref{ultra_rad3}.\dokendDef
\end{defin}

\begin{rem} {\bf Properties of weak entropy solutions}
\begin{itemize}
\item[1)] It follows from the assumptions in Definition \ref{def:weakrad}
for $a$, $b$, $a_0$ and $b_0$ that all the integrals in \eqref{ultra_rad2}
and \eqref{weak_entropy} are well defined. For \eqref{ultra_rad2} we first note that $x a$
is locally integrable. From $|b|<a$ and \eqref{cdefinition} we conclude that $x b$ 
and $x c(a,b)$ are locally integrable as well. For the entropy inequality
we make use of \eqref{statetrans} and obtain 
$$
p^{3/4}\sqrt{1+u^2} = 2^{-\frac{7}{4}}(3a-c(a,b))^{\frac12} (a-c(a,b))^{\frac14}\,.
$$
With $\frac{a}{3} \leq c(a,b) < a$ we have
$\begin{displaystyle}p^{3/4}|u|\leq p^{3/4}\sqrt{1+u^2}\end{displaystyle}\leq (a/3)^{3/4}$
and conclude that the integral in \eqref{weak_entropy} is well defined.
\item[2)] At all points of smoothness any solution $p,u$ satisfying \eqref{ultra_rad1}
will also satisfy the entropy conservation law
\begin{equation*}
\frac{\partial}{\partial t}\left(x^2 p^{3/4}\sqrt{1+u^2}\right) + \frac{\partial }{\partial x}
 \left(x^2 p^{3/4} u\right)= 0\,.
\end{equation*}
In this case the additional conservation law can be obtained from \eqref{ultra_rad1} by a straightforward but lengthy calculation. 
\item[3)] In the presence of shock waves \eqref{weak_entropy} will satisfy the strict inequality in general, 
and we obtain a simple evaluation of \eqref{weak_entropy}, see \cite[Chapter 2.1]{Abdelrahman1} for more details:
If for $p_-,p_+>0$ the left state $(p_-,u_-)$ can be connected to the right state $(p_+,u_+)$ by a single shock satisfying the Rankine-Hugoniot jump conditions, then this shock wave satisfies the entropy inequality if and only if $u_- > u_+$.
This condition can also be checked easily for the numerical solutions with shock curves in Chapter \ref{numeric}.
\end{itemize}
\end{rem}

In \cite{Kunikthesis}
we have used contour integrals for weak solutions of conservation laws,
following Oleinik's formulation \cite{Oleinik:1957} for a scalar conservation law.
Here we recall the definition \eqref{cdefinition} of $c$,
make use of the abbreviations $a=a(t,x)$, $b=b(t,x)$
and obtain an alternative formulation of \eqref{ultra_rad2} 
if we require especially for a piecewise smooth weak solution $a$, $b$ 
and for each convex domain $\Omega \subset Q$ 
with piecewise smooth boundary $\partial \Omega \subset Q$:
\begin{equation}\label{weak_contour}
\begin{array}{ll}
& \begin{displaystyle}
\int \limits_{\partial \Omega}x^2 a\,dx-x^2 b\,dt 
\end{displaystyle}= 0\,,\\
& \begin{displaystyle}
\int \limits_{\partial \Omega}x^2 b\,dx-x^2 c\,dt 
\end{displaystyle}= 
\begin{displaystyle}
\iint \limits_{\Omega} x(a-c)\,dtdx 
\end{displaystyle}\,.\\
\end{array}
\end{equation}
\section{Solutions of the linearized system}\label{linear}
A linearized version of the system \eqref{ultra_rad3} is given by
\begin{equation}\label{ultra_linearized}
\left\{
\begin{array}{ll}
 & \begin{displaystyle}\frac{\partial}{\partial t}\left(x^2 a \right)
+ \frac{\partial}{\partial x}\left(x^2 b \right)\end{displaystyle}=0\,,\\[1.1ex]
& \begin{displaystyle}\frac{\partial}{\partial t}\left(x^2 b \right)
+ \frac{\partial}{\partial x}
\left(\frac{x^2}{3} a \right)\end{displaystyle}
=\begin{displaystyle}\frac{2x}{3}a\,.\end{displaystyle}\\
\end{array}
\right.
\end{equation}
We linearize at the state $(a,b)=(a,0)$. 
The system can be obtained by neglecting the terms $b^2$ in \eqref{ultra_rad3}.
For $t=0$ and $x\geq 0$ we prescribe initial data $a_0(x)=a(0,x)$,
$b_0(x)=b(0,x)$ and assume that $b_0(0)=0$. 
From the radial symmetry the variable $x>0$ corresponds to the radius variable. 
Now we want to extend our initial data to all of $\R$ using symmetry in order to obtain simple solution formulas.
For $x>0$ we extend $a_0$ to an \textit{even function} with $a_0(-x)=a_0(x)$
and $b_0$ to an \textit{odd function} with $b_0(-x)=-b_0(x)\,.$
Now we assume that $a_0, b_0 : \R \mapsto \R$ are both $C^1$-functions. 
For $x \in \R$ we define the two \textit{even primitive functions}
\begin{equation}\label{primitives}
A_0(x)=\int \limits_0^x u a_0(u)\,du\,, \quad B_0(x)=\int \limits_0^x b_0(u)\,du\,.
\end{equation}

\begin{thm}\label{linsolv}
The solution of \eqref{ultra_linearized}
satisfying the initial conditions
$$
\lim \limits_{t \,\downarrow \, 0}a(t,x)=a_0(x)\,, \quad
\lim \limits_{t \,\downarrow \, 0}b(t,x)=b_0(x)\,, \quad x>0\,,
$$
is given for all $t>0$, $x>0$ by
\begin{equation}\label{solution_alpha}
\begin{split}
a(t,x)&=\\
&\frac{1}{2x}
\left(x+\frac{t}{\sqrt{3}} \right)a_0\left(x+\frac{t}{\sqrt{3}} \right)\\
+&\frac{1}{2x}
\left(x-\frac{t}{\sqrt{3}} \right)a_0\left(x-\frac{t}{\sqrt{3}} \right)\\
-&\frac{\sqrt{3}}{2x}
\left[
\left(x+\frac{t}{\sqrt{3}} \right)b_0\left(x+\frac{t}{\sqrt{3}} \right)+
B_0\left(x+\frac{t}{\sqrt{3}} \right) 
\right]\\
+&\frac{\sqrt{3}}{2x}\,\left[
\left(x-\frac{t}{\sqrt{3}} \right)b_0\left(x-\frac{t}{\sqrt{3}} \right)
+ B_0\left(x-\frac{t}{\sqrt{3}}\right)\right]\,,\\
\end{split}
\end{equation}
%%%%%%%%%%%%%%%%%%%%%%%%%%%%%%%%%%%%%%%%%%%%%%%%%%%%%%%%%%%
\begin{equation}\label{solution_beta}
\begin{split}
b(t,x)& =\\
-&\frac{1}{2\sqrt{3}x}
\left[
\left(x+\frac{t}{\sqrt{3}} \right)a_0\left(x+\frac{t}{\sqrt{3}} \right)
-\frac{1}{x}\,A_0\left(x+\frac{t}{\sqrt{3}} \right)
 \right]\\
+&\frac{1}{2\sqrt{3}x}\,\left[
\left(x-\frac{t}{\sqrt{3}} \right)a_0\left(x-\frac{t}{\sqrt{3}} \right)
- \frac{1}{x}\,A_0\left(x-\frac{t}{\sqrt{3}} \right)\right]\\
+&\frac{1}{2x}
\left[
\left(x+\frac{t}{\sqrt{3}} \right)b_0\left(x+\frac{t}{\sqrt{3}} \right)
-\frac{1}{x}\cdot\frac{t}{\sqrt{3}}\,B_0\left(x+\frac{t}{\sqrt{3}} \right) 
\right]\\
+&\frac{1}{2x}\,\left[
\left(x-\frac{t}{\sqrt{3}} \right)b_0\left(x-\frac{t}{\sqrt{3}} \right)
+\frac{1}{x}\cdot\frac{t}{\sqrt{3}}\,B_0\left(x-\frac{t}{\sqrt{3}}\right)\right]\,.\\
\end{split}
\end{equation}
\end{thm}
\begin{proof}
Assume that $F_0'=f_0$ for any $C^2$-function $F_0$.
If we make the ansatz
\begin{equation}\label{plus}
\begin{split}
a_+(t,x)&=\frac{1}{x}f_0\left(x+\frac{t}{\sqrt{3}} \right)\,,\\
b_+(t,x)&=-\frac{1}{\sqrt{3}\,x}\,f_0\left(x+\frac{t}{\sqrt{3}} \right)
+\frac{1}{\sqrt{3}\,x^2}\,F_0\left(x+\frac{t}{\sqrt{3}} \right)\,,\\
\end{split}
\end{equation}
then we can easily check that $a_+,b_+: Q \mapsto \R$ 
satisfy \eqref{ultra_linearized} with
$$
\lim \limits_{t \,\downarrow \, 0}a_+(t,x)=\frac{f_0(x)}{x}\,, \quad
\lim \limits_{t \,\downarrow \, 0}b_+(t,x)=-\frac{f_0(x)}{\sqrt{3}\,x}
+\frac{F_0(x)}{\sqrt{3}\,x^2}\,, \quad x>0\,.
$$
In the same way
\begin{equation}\label{minus}
\begin{split}
a_-(t,x)&=\frac{1}{x}f_0\left(x-\frac{t}{\sqrt{3}} \right)\,,\\
b_-(t,x)&=\frac{1}{\sqrt{3}\,x}\,f_0\left(x-\frac{t}{\sqrt{3}} \right)
-\frac{1}{\sqrt{3}\,x^2}\,F_0\left(x-\frac{t}{\sqrt{3}} \right)\,,\\
\end{split}
\end{equation}
also satisfy \eqref{ultra_linearized} for $t>0$, $x>0$ with
$$
\lim \limits_{t \,\downarrow \, 0}a_-(t,x)=\frac{f_0(x)}{x}\,, \quad
\lim \limits_{t \,\downarrow \, 0}b_-(t,x)=\frac{f_0(x)}{\sqrt{3}\,x}
-\frac{F_0(x)}{\sqrt{3}\,x^2}\,, \quad x>0\,.
$$
Using \eqref{primitives}, \eqref{plus} and \eqref{minus} we can check that \eqref{solution_alpha} 
and \eqref{solution_beta} satisfy \eqref{ultra_linearized} line by line
in the following way:

For the first line on the right-hand side in \eqref{solution_alpha} and 
\eqref{solution_beta} we put $F_0(x)=\frac12\,A_0(x)$ and use \eqref{plus},
for the second line we put again $F_0(x)=\frac12\,A_0(x)$ and use \eqref{minus},
and the third and fourth line with 
$F_0(x)=-\frac{\sqrt{3}}{2}\,x\,B_0(x)$, \eqref{plus}
and $F_0(x)=+\frac{\sqrt{3}}{2}\,x\,B_0(x)$, \eqref{minus}, respectively.
In this way the initial conditions are also obtained.
\end{proof}

\begin{rem}
It follows from the previous theorem that $x\,a(t,x)$ and $x\,b(t,x)$
are even bounded in any bounded subdomain of the quarterplane $t,x>0$.\\
\end{rem}

However, for more general weak solutions $|a(t,x)|$ and $|b(t,x)|$ 
may become infinitly large in certain small time intervals for $x \downarrow 0$. 
This is shown in the following example with a spherical imploding shock: 
We first use \eqref{primitives}.
Put $b_0(x)=B_0(x)=0$ for all $x \in \R$ and
\begin{equation*}
a_0(x)=\left\{
\begin{array}{ll}
1 & \text{~for~} |x| < 1\,,\\
2 & \text{~for~} |x| \geq 1\,.\\
\end{array}
\right.
\end{equation*}
Then we have the even function
\begin{equation*}
A_0(x)=\int \limits_{0}^{x} u a_0(u)du=\left\{
\begin{array}{ll}
\frac12 \, x^2 & \text{~for~} |x| < 1\,,\\
x^2-\frac12 & \text{~for~} |x| \geq 1\,.\\
\end{array}
\right.
\end{equation*}
We define, as seen in Figure \ref{fig1}, the convex domains
\begin{equation*}
\begin{split}
 \Omega_1&=\{\,(t,x) \in \R_{>0} \times \R_{>0}\,:\, x \leq 1-t/\sqrt{3}\}\,,\\
 \Omega_2&=\{\,(t,x) \in \R_{>0} \times \R_{>0}\,:\, |1-t/\sqrt{3}|< x \leq 1+t/\sqrt{3}\}\,,\\
 \Omega_3&=\{\,(t,x) \in \R_{>0} \times \R_{>0}\,:\, x>  1+t/\sqrt{3}\}\,,\\
 \Omega_4&=\{\,(t,x) \in \R_{>0} \times \R_{>0}\,:\, x \leq  t/\sqrt{3}-1\}\\
\end{split}
\end{equation*}
and obtain from \eqref{solution_alpha}, \eqref{solution_beta}:
\begin{equation*}
a(t,x)= \left\{\begin{array}{ccl}
 1   &\quad  &\mbox{~for~} (t,x) \in \Omega_1\,,\\[1ex]
\begin{displaystyle}\frac{3}{2}+\frac{t}{2\sqrt{3}\,x}
\end{displaystyle}&\quad  &\mbox{~for~} (t,x) \in \Omega_2\,,\\[2ex]
2    &\quad &\mbox{~for~} (t,x) \in \Omega_3 \cup \Omega_4\,,\\
\end{array}
\right.
\end{equation*}
\begin{equation*}
b(t,x)= \left\{\begin{array}{ccl}
 0   &\quad  &\mbox{~for~} (t,x) \in \Omega_1 \cup \Omega_3 \cup \Omega_4\,,\\[1ex]
\begin{displaystyle}\frac{t^2-3(1+x^2)}{12\sqrt{3}\,x^2} 
\end{displaystyle}&\quad  &\mbox{~for~} (t,x) \in \Omega_2\,.\\
\end{array}
\right.
\end{equation*}

\begin{figure}
\scalebox{0.87}{%\documentclass{article}
%\usepackage{rotating}
%\usepackage{pst-all}
%\begin{document}
%\pagestyle{empty}
%\begin{center}
\begin{pspicture}(0,0)(10,10)

\psline{-*}(0,0)(0,0)
\rput[t]{0}(-0.3,-0.3){$0$}
\psline{-*}(0,3)(0,3)
\rput[t]{0}(-0.3,3){$1$}
\psline[linewidth=1pt]{->}(0,0)(11,0)
\rput[t]{0}(11,-0.3){$t$}
\psline[linewidth=1pt]{->}(0,0)(0,6.5)
\rput[t]{0}(-0.3,6.5){$x$}

\psline{-*}(0,3)(5.2,0)
\psline(5.2,0)(10.4,3)
\psline(0,3)(5.2,6)

\pcline[linestyle=none](0,3)(5.2,0)
\aput{:U}{\mbox{$x=x_1(t)$}}
\pcline[linestyle=none](5.2,0)(10.4,3)
\aput{:U}{\mbox{$x=x_2(t)$}}
\pcline[linestyle=none](0,3)(5.2,6)
\aput{:U}{\mbox{$x=x_3(t)$}}

\rput[t]{0}(5.2,-0.1){$\sqrt{3}$}
\rput[t]{0}(0.8,0.6){$\Omega_1$}
\rput[t]{0}(1.5,1.3){$a=1, b=0$}

\rput[t]{0}(1.6,3.0){$\Omega_2$}
\rput[t]{0}(4.7,4.5){$\begin{displaystyle}a=\frac{3}{2}+\frac{t}{2\sqrt{3}\,x}\,,
\end{displaystyle}$}
\rput[t]{0}(4.9,3.3){$\begin{displaystyle}b=\frac{t^2-3(1+x^2)}{12\sqrt{3}\,x^2}
\end{displaystyle}$}

\rput[t]{0}(0.8,4.5){$\Omega_3$}
\rput[t]{0}(1.5,5.9){$a=2, b=0$}

\rput[t]{0}(7.2,0.6){$\Omega_4$}
\rput[t]{0}(9.2,1.3){$a=2, b=0$}
\end{pspicture}
%\end{center}
%\end{document}
}
\caption{Solution of the linearized system.}
\label{fig1}
\end{figure}

We can easily check that $a, b$ satisfy the differential equations
\eqref{ultra_linearized} in the interior of each $\Omega_j$, $j=1,2,3,4$.
Moreover, for $i=1,2,3$ the Rankine-Hugoniot jump conditions 
$$
\dot{x_i}(t)=\lim \limits_{\varepsilon \,\downarrow \, 0}
\frac{b(t,x_i(t)+\varepsilon)-b(t,x_i(t)-\varepsilon)}
{a(t,x_i(t)+\varepsilon)-a(t,x_i(t)-\varepsilon)}
=\pm \frac{1}{\sqrt{3}}
$$
of the linearized system are satisfied across the three shocks
\begin{align*}
x_1(t)&=1-\frac{t}{\sqrt{3}}\,, \quad 0<t<\sqrt{3}\,,\\
x_2(t)&=\frac{t}{\sqrt{3}}-1\,, \quad t>\sqrt{3}\,,\\
x_3(t)&=1+\frac{t}{\sqrt{3}}\,, \quad t>0\,,\\
\end{align*}
see Figure \ref{fig1}. 
Theorem \ref{linsolv} gives a weak solution of the linearized system
even if the initial functions $a_0$, $b_0$ have jump discontinuities.
\begin{thm}\label{linboundaries}
Take the assumptions as in Theorem \ref{linsolv}. If $b_0, a_0$ are both $C^2$-functions, then we have \,
$\begin{displaystyle}
%&\lim \limits_{x\,\downarrow\,0}a(t,x)=
%a_0(t/\sqrt{3})-2\sqrt{3}b_0(t/\sqrt{3})+t/\sqrt{3}\,a'_0(t/\sqrt{3})
%-tb_0'(t/\sqrt{3})\,,\\
\lim \limits_{x\,\downarrow\,0}b(t,x)=0 \,
\end{displaystyle}$
for all $t>0$.
\end{thm}
\begin{proof}
L'Hospital's rule can be applied twice to obtain the desired radial limit.
\end{proof}
The linearized system serves as a motivation for the following study
of the nonlinear system. However, we cannot expect a quantitativly similar behaviour between both models 
concerning the shock-wave reflection and the singular structure near the boundary,
where nonlinear momentum terms cannot be neglected.

%%%%%%%%%%%%%%%%%%%%%%%%

\section{Formulation of a stable numerical scheme}\label{scheme1}
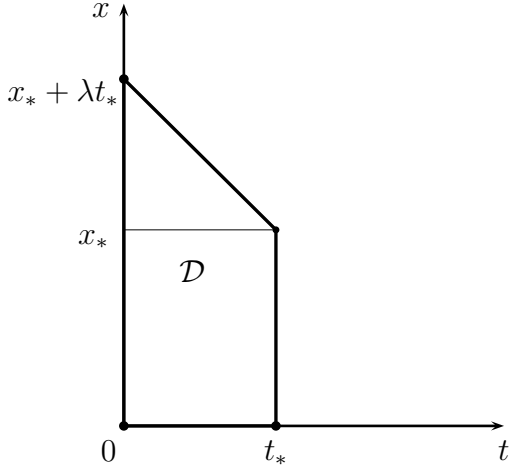
\begin{figure}
\scalebox{1}{\begin{pspicture}(0,0)(6,6)

\psline{-*}(0,0)(0,0)
\rput[t]{0}(-0.2,-0.2){$0$}
\rput[t]{0}(-0.4,2.6){$x_*$}
\psline{-*}(0,4.6)(0,4.6)
\psline{-*}(2,0)(2,0)

\rput[t]{0}(-0.8,4.6){$x_{*}+\lambda t_*$}
\rput[t]{0}(0.9,2.2){$\mathcal{D}$}
\psline[linewidth=1pt]{->}(0,0)(5,0)
\rput[t]{0}(5,-0.2){$t$}

\psline[linewidth=1pt]{->}(0,0)(0,5.6)
\rput[t]{0}(-0.3,5.6){$x$}

\psline[linewidth=1.25pt]{-}(0,0)(2,0)
\psline[linewidth=1.25pt]{-}(2,0)(2,2.6)
\psline[linewidth=1.25pt]{-}(2,2.6)(0,4.6)
\psline[linewidth=1.25pt]{-}(0,0)(0,4.6)
\rput[t]{0}(2,-0.2){$t_*$}

\psline[linewidth=0.3pt]{-*}(0,2.6)(2,2.6)
\psline(0,4.6)(2,2.6)
\end{pspicture}}
\caption{The computational domain $\mathcal{D}$}
\label{fig_domain}
\end{figure}

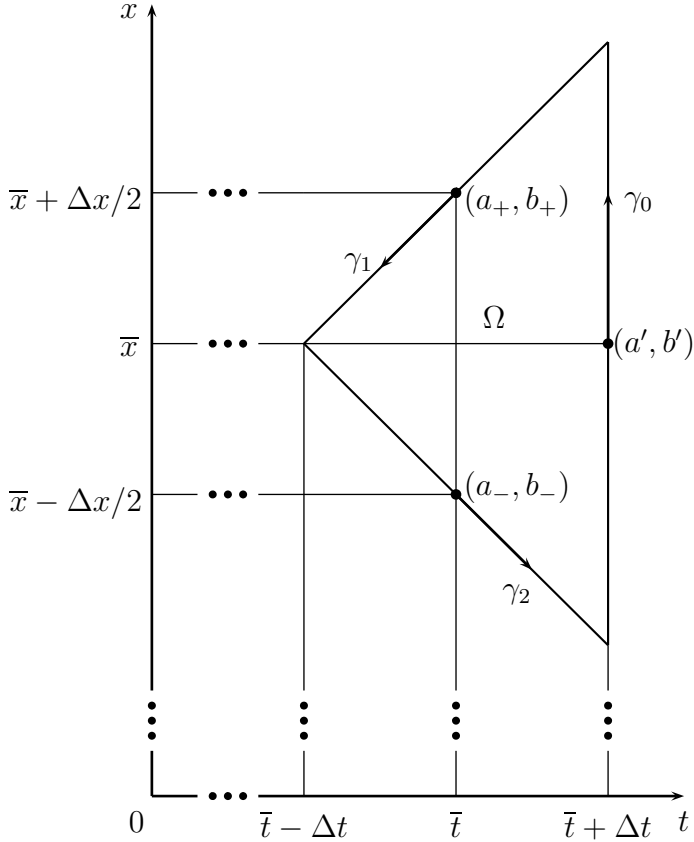
\begin{figure}
\scalebox{1}{\begin{pspicture}(-2.5,-2.5)(9,9)

\psline[linewidth=0.5pt]{-}(0,2)(0.6,2)
\psline[linewidth=0.5pt]{-*}(0.8,2)(0.8,2)
\psline[linewidth=0.5pt]{-*}(1,2)(1,2)
\psline[linewidth=0.5pt]{-*}(1.2,2)(1.2,2)
\psline[linewidth=0.5pt]{-}(1.4,2)(2,2)

\psline[linewidth=0.5pt]{-}(0,4)(0.6,4)
\psline[linewidth=0.5pt]{-*}(0.8,4)(0.8,4)
\psline[linewidth=0.5pt]{-*}(1,4)(1,4)
\psline[linewidth=0.5pt]{-*}(1.2,4)(1.2,4)
\psline[linewidth=0.5pt]{-}(1.4,4)(2,4)

\psline[linewidth=0.5pt]{-}(0,6)(0.6,6)
\psline[linewidth=0.5pt]{-*}(0.8,6)(0.8,6)
\psline[linewidth=0.5pt]{-*}(1,6)(1,6)
\psline[linewidth=0.5pt]{-*}(1.2,6)(1.2,6)
\psline[linewidth=0.5pt]{-}(1.4,6)(2,6)

\rput[t]{0}(-0.2,-2.2){$0$}
\rput[t]{0}(-1,2.1){$\overline{x}-\Delta x/2$}
\psline{*-*}(4,6)(4,6)
\rput[t]{0}(-1,6.1){$\overline{x}+\Delta x/2$}
\rput[t]{0}(-0.3,4.1){$\overline{x}$}
%\rput[t]{0}(4.6,5){$\Omega_{n,j}^+$}
%\rput[t]{0}(3.4,5){$\Omega_{n,j}^-$}

%\rput[t]{0}(3.8,4.3){$a=constant$}
%\rput[t]{0}(3.8,3.9){$b=constant$}
%\rput[t]{0}(2.75,2.75){$a_-,b_-$}
%\rput[t]{0}(2.75,5.75){$a_+,b_+$}

\psline[linewidth=1pt]{-}(0,-2)(0.6,-2)
\psline[linewidth=0.5pt]{-*}(0.8,-2)(0.8,-2)
\psline[linewidth=0.5pt]{-*}(1,-2)(1,-2)
\psline[linewidth=0.5pt]{-*}(1.2,-2)(1.2,-2)
\psline[linewidth=1pt]{->}(1.4,-2)(7,-2)
\rput[t]{0}(7,-2.2){$t$}

\psline[linewidth=1pt]{-}(0,-2)(0,-1.4)
\psline[linewidth=0.5pt]{-*}(0,-1.2)(0,-1.2)
\psline[linewidth=0.5pt]{-*}(0,-1)(0,-1)
\psline[linewidth=0.5pt]{-*}(0,-0.8)(0,-0.8)
\psline[linewidth=1pt]{->}(0,-0.6)(0,8.5)
\rput[t]{0}(-0.3,8.5){$x$}

\psline[linewidth=0.5pt]{-}(2,-2)(2,-1.4)
\psline[linewidth=0.5pt]{-*}(2,-1.2)(2,-1.2)
\psline[linewidth=0.5pt]{-*}(2,-1)(2,-1)
\psline[linewidth=0.5pt]{-*}(2,-0.8)(2,-0.8)
\psline[linewidth=0.5pt]{-}(2,-0.6)(2,4)

\psline[linewidth=0.5pt]{-}(4,-2)(4,-1.4)
\psline[linewidth=0.5pt]{-*}(4,-1.2)(4,-1.2)
\psline[linewidth=0.5pt]{-*}(4,-1)(4,-1)
\psline[linewidth=0.5pt]{-*}(4,-0.8)(4,-0.8)
\psline[linewidth=0.5pt]{-}(4,-0.6)(4,2)

\psline[linewidth=0.5pt]{-}(6,-2)(6,-1.4)
\psline[linewidth=0.5pt]{-*}(6,-1.2)(6,-1.2)
\psline[linewidth=0.5pt]{-*}(6,-1)(6,-1)
\psline[linewidth=0.5pt]{-*}(6,-0.8)(6,-0.8)
\psline[linewidth=0.5pt]{-}(6,-0.6)(6,0)

\psline[linewidth=0.5pt]{-*}(2,4)(6,4)

\psline[linewidth=0.5pt]{-*}(2,2)(4,2)
\psline[linewidth=0.5pt]{-}(2,6)(4,6)
\psline[linewidth=0.5pt]{-}(4,2)(4,6)

\rput[t]{0}(2,-2.2){$\overline{t}-\Delta t$}
\rput[t]{0}(4,-2.2){$\overline{t}$}

\rput[t]{0}(6,-2.2){$\overline{t}+\Delta t$}

\psline{-}(2,4)(6,8)
\psline{-}(2,4)(6,0)
\psline{-}(6,0)(6,8)

\psline[linewidth=1pt]{*->}(4,2)(5,1)
\rput[t]{0}(4.8,0.8){$\gamma_2$}
\psline[linewidth=1pt]{*->}(4,6)(3,5)
\rput[t]{0}(2.72,5.2){$\gamma_1$}

\psline[linewidth=1pt]{*->}(6,4)(6,6)
\rput[t]{0}(6.4,6){$\gamma_0$}

\rput[t]{0}(4.5,4.5){$\Omega$}
\rput[t]{0}(4.8,6.1){$(a_+,b_+)$}
\rput[t]{0}(4.8,2.3){$(a_-,b_-)$}
\rput[t]{0}(6.6,4.2){$(a',b')$}

\end{pspicture}}
\caption{The balance region $\Omega$}
\label{figeul}
\end{figure}

We develop a stable numerical scheme 
for the initial value problem with the radially symmetric 
ultra-relativistic Euler equations.
The method of contour-integration for the formulation 
of the balance laws \eqref{weak_contour} is used to construct 
a function called "Euler" which enables the evolution in time 
of the numerical solution on a staggered grid, i.e. it allows us 
to construct the solution $(a',b')$ at the next time step from the solution $(a_{\pm},b_{\pm})$
in two neighboring gridpoints at the former time step according to Figure \ref{figeul}.
First we determine the computational domain and define some quantities 
which are needed for its discretization.
\begin{itemize}
\item[1)] Given are $t_*,x_*>0$ in order to calculate a numerical solution 
of the initial value problem \eqref{weak_contour}, \eqref{initial_final} in the time range $[0,t_*]$ 
and the spatial range $[0,x_*]$.
\item[2)] We want to use a staggered grid scheme. Any given number $N \in \N$ with $N \cdot x_* \geq t_*$ determines the time step size
\begin{equation*}
\Delta t = \frac{t_*}{2N}\,.
\end{equation*}
The time steps are
\begin{equation*}
t_{n} = (n-1) \Delta t\,, \quad n=1,\ldots,2N+1\,.
\end{equation*}
\item[3)] Put
\begin{equation*}
M = \left \lfloor \frac{x_*}{t_*}\,N \right \rfloor \geq 1\,,
\end{equation*}
then the spatial mesh size is
\begin{equation*}
\Delta x=\frac{x_*}{M}\,,
\end{equation*}
with the spatial grid points
\begin{equation*}
x_{j} = (j-1) \Delta x\,, \quad j=1,\ldots,N+M+1\,.
\end{equation*}
Note that our scheme uses a trapezoidal computational domain $\mathcal{D}$
defined below that includes the target domain $[0,t_*] \times [0,x_*]$.
Thereby, we can use all initial data that influence the solution on the
target domain. In this way we avoid using a numerical boundary condition at $x_*$.

\item[4)] The number
\begin{equation*}
 \lambda=\frac{\Delta x}{2 \Delta t} \geq 1
\end{equation*}
is used to satisfy the CFL-condition and to define the computational domain
$\begin{displaystyle}
 \mathcal{D} = \left \{(t,x)\in \R^2\,:\,0 \leq t \leq t_* \,,\quad 
0 \leq x \leq x_*+\lambda(t_*-t)\,\right \}\,.
\end{displaystyle}$
\end{itemize}
The typical trapezoidal form of the computational domain is illustrated in Figure \ref{fig_domain}.

For the formulation and the stability of our numerical scheme we need two lemmas.
\begin{lem}\label{hilfssatz1}
Assume that $|b_{\pm}|< a_{\pm}$. We recall that $\lambda\geq 1$,
\eqref{cdefinition} and put $c_{\pm}=c(a_{\pm},b_{\pm})$\,. Then
\begin{itemize}
\item[a)] $\begin{displaystyle}
-\left(a_-+\frac{b_-}{\lambda} \right)< b_-+\frac{c_-}{\lambda} < a_-+\frac{b_-}{\lambda}\,,
\end{displaystyle}$
\item[b)] $\begin{displaystyle}
-\left(a_+-\frac{b_+}{\lambda} \right)<  b_+-\frac{c_+}{\lambda} < a_+-\frac{b_+}{\lambda}\,.
\end{displaystyle}$
\end{itemize}
\end{lem}
\begin{proof}
The proof of b) is quite analogous to a), hence we will only show a).
The left inequality in a) is equivalent with
$\begin{displaystyle}
\lambda(a_- +b_-)>-c_- - b_-\,,
\end{displaystyle}$
and due to $\lambda\geq1$ and $a_- +b_->0$ it is sufficient to show
$\begin{displaystyle}
a_- +b_->-c_- - b_-\,.
\end{displaystyle}$
We have
\begin{align*}
& \qquad a_- +b_->-c_- - b_-\\
&\Leftrightarrow 2a_-+\frac32 b_->\sqrt{a_-^2-\frac34b_-^2} \quad \mbox{with~}2a_-+\frac32 b_->0\\
&\Leftrightarrow 4a_-^2 +6a_-b_- +\frac{9}{4}b_-^2>a_-^2-\frac34b_-^2\\
&\Leftrightarrow 3(a_-+b_-)^2>0\,,\\
\end{align*}
which shows the left inequality.
The right inequality in a) is equivalent to
$\begin{displaystyle}
\lambda(a_- -b_-)>c_- - b_-\,.
\end{displaystyle}$
Due to $\lambda\geq1$ and $a_- -b_->0$ it is sufficient to show
$\begin{displaystyle}
a_->c_-\,.
\end{displaystyle}$
We have
\begin{align*}
& \qquad a_->c_- \\
&\Leftrightarrow 2\sqrt{a_-^2-\frac34b_-^2}>a_-\quad \mbox{with~}a_->0\\
&\Leftrightarrow 4a_-^2 -3b_-^2>a_-^2\\
&\Leftrightarrow 3(a_-^2-b_-^2)>0\,.\\
\end{align*}
\end{proof}

\begin{lem}\label{hilfssatz2}
Assume that $a>0$, $0 < \eta \leq 1/3$ and $-a(1+\eta)<\xi<a(1-\eta)$. 
Then we obtain $4a^2(1+3\eta^2)-3\xi^2>0$ and
\begin{equation*}
\left|\frac{\xi + \eta\sqrt{4a^2(1+3\eta^2)-3\xi^2}}{1+3\eta^2}\right|<a\,.
\end{equation*}
\end{lem}
\begin{proof} We have
\begin{align*}
4a^2(1+3\eta^2)-3\xi^2>4a^2(1+3\eta^2)-3a^2(1+\eta)^2=a^2(1-3\eta)^2\geq 0\,,
\end{align*}
and the square root in the estimate of the lemma is well defined. To show the estimate we 
use $0 < \eta \leq 1/3$ and first note that
\begin{equation*}
\begin{split}
&\frac{-a(1+\eta) + \eta\sqrt{4a^2(1+3\eta^2)-3(-a(1+\eta))^2}}{1+3\eta^2}=-a\,,\\
&\frac{a(1-\eta) + \eta\sqrt{4a^2(1+3\eta^2)-3(a(1-\eta))^2}}{1+3\eta^2}=a\,.\\
\end{split}
\end{equation*}
Therefore it is sufficient for the proof of the lemma to show that
$$
B(\xi)=\xi + \eta\sqrt{4a^2(1+3\eta^2)-3\xi^2}
$$
is strictly monotonically increasing for $-a(1+\eta)<\xi<a(1-\eta)$.
The condition $B'(\xi)=0$ gives
$\begin{displaystyle}
3\eta\xi=\sqrt{4a^2(1+3\eta^2)-3\xi^2}\,,
\end{displaystyle}$
hence $\xi>0$ and the unique solution $\xi=\frac{2a}{\sqrt{3}}>a(1-\eta)$ outside the interval.
On the other hand we have $B'(0)=1>0$. Hence $B$ is strictly increasing in the interval.
\end{proof}

For the numerical discretization of the system \eqref{weak_contour}
we choose the triangular balance domain $\Omega$ depicted in Figure \ref{figeul}.
We assume that the midpoints $P_-=(\overline{t}, \overline{x}-\Delta x/2)$,
$P_+=(\overline{t}, \overline{x}+\Delta x/2)$ and 
$P'=(\overline{t}+\Delta t, \overline{x})$ of the cords of $\partial \Omega$
are numerical gridpoints for the computational domain $\mathcal{D}$. 
Let the numerical solution $(a_\pm,b_\pm)$
be given at the gridpoints $P_{\pm}$. 
We have to require $|b_{\pm}|< a_{\pm}$ for the numerical solution
in the actual time step $\overline{t}=t_n$ with $t=1,\ldots,2N$. 
The major task is to calculate
the numerical solution $(a',b')$ for the next time step 
$\overline{t}+\Delta t= t_{n+1}$ at its gridpoint $P'$, see Figure \ref{figeul}.

The spatial value $\overline{x} \geq 0$ is given.
We have to determine a function
\begin{equation}\label{euler_solve}
\text{Euler}(a_-,b_-,a_+,b_+,\overline{x},\Delta x, \lambda)=(a',b')
\end{equation}
for the calculation of $(a',b')$. This leads to the structure
of a staggered grid scheme. Note that at the boundary the balance region $\Omega$
may have parts outside $\mathcal{D}$, e.g. points
below the half-space $x\geq 0$. In the latter case we will employ a simple
reflection principle for the numerical solution in order to use the function
$\text{Euler}$ as well for the evaluation of the boundary conditions.

Next we will make use of the fact that the points
$P_{\pm}$ with numerical values $(a_{\pm},b_{\pm})$ and $P'$ 
with unknown value $(a',b')$ are the 
\textit{midpoints of the three boundary cords}
of the balance region $\Omega$. 
We put $c_{\pm}=c(a_{\pm},b_{\pm})$ and $c'=c(a',b')$
for abbreviation, see \eqref{cdefinition}.
Then we use for $k=0,1,2$ the straight line paths $\gamma_k$ 
from Figure \ref{figeul} and for the corresponding path integrals
$$
\int\limits_{\gamma_k}x^2 a(t,x)\,dx-x^2 b(t,x)\,dt \mbox{~and ~} 
\int\limits_{\gamma_k}x^2 b(t,x)\,dx-x^2 c(a(t,x),b(t,x))\,dt 
$$
with the unknown weak entropy solution $a(t,x)$, $b(t,x)$
their numerical discretizations $I_{k,a}$ and $I_{k,b}$, respectively, given by
\begin{equation}\label{pathint0}
\begin{split}
I_{0,a}&=\int\limits_{\gamma_0}x^2 a'\,dx-x^2 b'\,dt=a'\,
\int \limits_{\overline{x}-\Delta x}^{\overline{x}+\Delta x}x^2\,dx=
4a'\lambda \Delta t \left( \overline{x}^2 +\frac13 (\Delta x)^2 \right)\,,\\
I_{0,b}&=\int\limits_{\gamma_0}x^2 b'\,dx-x^2 c'\,dt=b'\,
\int \limits_{\overline{x}-\Delta x}^{\overline{x}+\Delta x}x^2\,dx=
4b'\lambda \Delta t \left( \overline{x}^2 +\frac13 (\Delta x)^2 \right)\,.\\
\end{split}
\end{equation}
\begin{equation}\label{pathint1}
\begin{split}
I_{1,a}&=\int\limits_{\gamma_1}x^2 a_+\,dx-x^2 b_+\,dt=
-2(\lambda a_+-b_+) \Delta t 
\left( \overline{x}^2 +\frac13 (\Delta x)^2 +\overline{x}\Delta x\right)\,,\\
I_{1,b}&=\int\limits_{\gamma_1}x^2 b_+\,dx-x^2 c_+\,dt=
-2(\lambda b_+-c_+) \Delta t 
\left( \overline{x}^2 +\frac13 (\Delta x)^2 +\overline{x}\Delta x\right)\,.\\
\end{split}
\end{equation}
\begin{equation}\label{pathint2}
\begin{split}
I_{2,a}&=\int\limits_{\gamma_2}x^2 a_-\,dx-x^2 b_-\,dt=
-2(\lambda a_- +b_-) \Delta t 
\left( \overline{x}^2 +\frac13 (\Delta x)^2 -\overline{x}\Delta x\right)\,,\\
I_{2,b}&=\int\limits_{\gamma_2}x^2 b_-\,dx-x^2 c_-\,dt=
-2(\lambda b_- +c_-) \Delta t 
\left( \overline{x}^2 +\frac13 (\Delta x)^2 -\overline{x}\Delta x\right)\,.\\
\end{split}
\end{equation}
We recall that $\overline{x}\geq 0$ and put
\begin{equation}\label{qdef}
q = \frac{2 \overline{x}\,\Delta x}{\overline{x}\,^2+\frac{1}{3}(\Delta x)^2}<2\,.
\end{equation}
The numerical discretization of the first balance law in \eqref{weak_contour} gives
\begin{equation}\label{abilanz}
I_{0,a}=-I_{1,a}-I_{2,a}\,.
\end{equation}
We obtain from \eqref{pathint0}, \eqref{pathint1},\eqref{pathint2},
\eqref{qdef} and \eqref{abilanz} for $a'$ the explicit solution
\begin{equation}\label{aC}
a'=\frac12\left(a_- +\frac{b_-}{\lambda}\right)(1-q/2)+
\frac12\left(a_+ - \frac{b_+}{\lambda}\right)(1+q/2)\,.
\end{equation}
For the numerical discretization of the second balance law in \eqref{weak_contour}
we approximate the integral
\begin{equation*}
\iint \limits_{\Omega}x(a-c)dtdx
\end{equation*}
by 
\begin{equation*}
(a'- c')\iint \limits_{\Omega}xdtdx
=2(a'- c')\overline{x}\Delta t \Delta x\,.
\end{equation*}
Now \eqref{pathint0}, \eqref{pathint1}, \eqref{pathint2}
give the following ansatz for the calculation of $b'$:
\begin{equation}\label{bbilanz}
I_{0,b}=-I_{1,b}-I_{2,b}+2(a'- c')\overline{x}\Delta t \Delta x\,.
\end{equation}
Recall $c_\pm=c(a_\pm,b_\pm)$ with $c=c(a,b)$ in \eqref{cdefinition}.
We use the abbreviations
\begin{equation}\label{Cxieta}
\xi=\frac12\left(b_- +\frac{c_-}{\lambda}\right)(1-q/2)+
\frac12\left(b_+ - \frac{c_+}{\lambda}\right)(1+q/2)-\frac{a'q}{6\lambda}\,, \quad
\eta=\frac{q}{6\lambda}\,.
\end{equation}
From \eqref{bbilanz} we obtain the implicit equation
\begin{equation}\label{bCimp}
b'=\xi+\eta \sqrt{4a'^2-3b'^2}\,.
\end{equation}
This leads to a quadratic equation for $b'$.
Lemma \ref{hilfssatz1} gives 
$$-a'(1+\eta)<\xi<a'(1-\eta)$$
for the quantity $a'$ in \eqref{aC}.
In order to apply Lemma \ref{hilfssatz2} with $a'$ instead of $a$
we have to choose the solution
\begin{equation}\label{bCexplicit}
b'=\frac{\xi + \eta\sqrt{4a'^2(1+3\eta^2)-3\xi^2}}{1+3\eta^2}
\end{equation}
of \eqref{bCimp} with the positive square root.
Now $b'$ is well defined with $|b'|<a'$, 
see the transformation \eqref{statetrans} in state space. 
We summarize our results in the following

\begin{thm}\label{apbp} 
{\bf Numerical solution $(a',b')$ for the balance region $\Omega$}\\
Given are real quantities $\overline{x}\geq 0$ and $a_{\pm}$, $b_{\pm}$.
Assume that $|b_{\pm}|< a_{\pm}$. We recall 
$\lambda\geq 1$ defined in terms of $\Delta t$ and $\Delta x$
and put $c_{\pm}=c(a_{\pm},b_{\pm})$ in the definition \eqref{cdefinition}\,.
Then we have $|b'|<a'$ for the quantities $a'$ and $b'$ 
calculated from \eqref{qdef}, \eqref{aC}, \eqref{Cxieta} and \eqref{bCexplicit}.
\dokendDef
\end{thm}

\begin{defin}\label{eulerfunction} {\bf The function Euler}\\
The state $(a',b')$ from Theorem \ref{apbp}
defines the function $\text{Euler}$ in \eqref{euler_solve}.
\dokendDef
\end{defin}

\begin{rem}\label{boundrem}
Assume that the state $(a_+,b_+)$ with $|b_+|<a_+$ is given and
that $\overline{x}=0$. We define 
the "reflected state" ${\mathcal R}(a_+,b_+)=(a_+,-b_+)$ and obtain 
\begin{equation}\label{eulreflect}
\text{Euler}(a_+,-b_+,a_+,b_+,0,\Delta x,\lambda)=(a_+-b_+/\lambda,\,0)\,.\\
\end{equation}
This means that numerical values $(a',b')$ calculated
with the function $\text{Euler}$ in \eqref{euler_solve}
at the boundary $\overline{x}=0$ with reflected states 
$(a_-,b_-)={\mathcal R}(a_+,b_+)$
satisfy the boundary condition $b'=0$.
\end{rem}

Now we are able to formulate the numerical scheme for the solution of the initial-boundary value problem \eqref{initial_final}, \eqref{boundary_final}, \eqref{weak_contour}. We  construct staggered grid points in the computational domain 
$\mathcal{D}$ and compute the numerical solution at these gridpoints. 
The function $\text{Euler}$ enables the evolution 
of the numerical solution in time, i.e. it allows us
to construct the solution at time $t=t_{n+1}$ 
from the solution which is already calculated in the gridpoints 
at the former time step $t=t_n$.
Note that the triangular balance domains that we used 
to determine the routine $\text{Euler}$ overlap. But this presents no problem since they are not needed once the formulas for new values have been obtained.
\begin{itemize}
\item The staggered gridpoints are 
$(t_n,x_{n,j}) \in \mathcal{D}$ for $t_n=(n-1)\Delta t$, 
$n=1,\ldots,2N+1$ and $j=1,\ldots,M+N-\lfloor (n-1)/2 \rfloor$ with
$$
x_{n,j}=\left\{
\begin{array}{ll}
(x_{j}+x_{j+1})/2 & \mbox{if $n$ is odd}\\
x_j   & \mbox{if $n$ is even}\,.\\
\end{array}
\right.
$$
We want to calculate the numerical solution 
$(a_{n,j},b_{n,j})$ at $(t_n,x_{n,j})$.
\item For $j=1,\ldots,M+N$ we calculate the numerical solution
$(a_{1,j}, b_{1,j})$ at the gridpoint $(t_1,x_{1,j})=(0,(x_j+x_{j+1})/2)$ 
from the given initial data by
$$a_{1,j}=a_0\left(x_{1,j}\right)\,, \quad 
b_{1,j}=b_0\left(x_{1,j}\right)\,.$$
This corresponds to taking the integral average
of the initial data on $(x_j,x_{j+1})$ and using the midpoint rule
as quadrature.
\item Assume that for a fixed \textit{odd index} $n\in \{1,\ldots,2N\}$ we have already determined the numerical solution $(a_{n,j}, b_{n,j})$ 
at the gridpoints $(t_n,x_{n,j})$, $j=1,\ldots,M+N-(n-1)/2$.\\
First we determine the solution $(a_{n+1,1},b_{n+1,1})$ at the boundary point
$(t_{n+1},x_{n+1,1})=(t_{n+1},0)$ according to \eqref{eulreflect} in
Remark \ref{boundrem}.
For this purpose we put $a_+=a_{n,1}$, $b_+=b_{n,1}$, $a_-=a_{n,1}$, 
$b_-=-b_{n,1}$ and have
$$
(a_{n+1,1}, b_{n+1,1})=\text{Euler}(a_-,b_-,a_+,b_+,0,\Delta x,\lambda)
\mbox{~with~} b_{n+1,1}=0\,.
$$
Next we put $a_-=a_{n,j-1}$, $b_-=b_{n,j-1}$ and $a_+=a_{n,j}$, $b_+=b_{n,j}$ for 
$j=2,\ldots,M+N-(n-1)/2$ and determine the values 
$a_{n+1,j}$, $b_{n+1,j}$ 
at time $t_{n+1}$ and position $\overline{x}=x_{n+1,j}=x_{j}$ from
$$
(a_{n+1,j}, b_{n+1,j})=\text{Euler}(a_-,b_-,a_+,b_+,\overline{x},\Delta x,\lambda)\,.
$$
\item Assume that for a fixed \textit{even index} $n\in \{1,\ldots,2N\}$ we have already determined the numerical solution $(a_{n,j}, b_{n,j})$ 
at the gridpoints $(t_n,x_{n,j})$, $j=1,\ldots,M+N-n/2+1$.\\
We put $a_-=a_{n,j}$, $b_-=b_{n,j}$ and $a_+=a_{n,j+1}$, $b_+=b_{n,j+1}$ for 
$j=1,\ldots,M+N-n/2$ and determine the values 
$a_{n+1,j}$, $ b_{n+1,j}$ at time $t_{n+1}$ and position 
$\overline{x}=x_{n+1,j}=(x_{j}+x_{j+1})/2$ from
$$
(a_{n+1,j}, b_{n+1,j})=\text{Euler}(a_-,b_-,a_+,b_+,\overline{x},\Delta x,\lambda)\,.
$$
\end{itemize}
Based on Lemma \ref{hilfssatz1} and \ref{hilfssatz2} 
we obtained Theorem \ref{apbp}. This implies stability for our scheme,
namely the following
\begin{thm}\label{stable_scheme} 
The numerical scheme described above is stable, especially the numerical values for the pressure $p$ always remain  positive.\dokendDef
\end{thm}

\section{Numerical examples}\label{numeric}
We solve the initial value problem \eqref{ultra_rad3}, \eqref{initial_final}
numerically for different choices of the initial data $a_0,b_0$.
We make use of the transformation \eqref{statetrans}.
However, for our numerical results we take the usual velocity 
$$ v=\frac{u}{\sqrt{1+u^2}} \quad \mbox{with~}|v|<1$$
instead of the four velocity $u$ and the initial velocity $v_0=v(0,\cdot)$.
The restriction $|v|<1$ leads to better color plots.
\begin{itemize}
\item[1)] If $a_0>0$ is constant and $b_0=0$, then we obtain a stationary solution, 
which is exactly reconstructed with these values by the scheme 
in Section \ref{scheme1}. This corresponds to $u=v=0$ and constant pressure $p$.
Such a steady part is contained in the following examples.
\item[2)] We choose the constant initial data $a_0=7$, $b_0=4\sqrt{2}$ corresponding to a constant initial pressure $p_0=1$ and a constant radial part 
$u_0=1$ and $v_0=1/\sqrt{2}$ of the initial four velo\-city and usual velo\-city,
respectively.
\begin{center}
\begin{figure}[ht]
\includegraphics[width=0.75\textwidth]{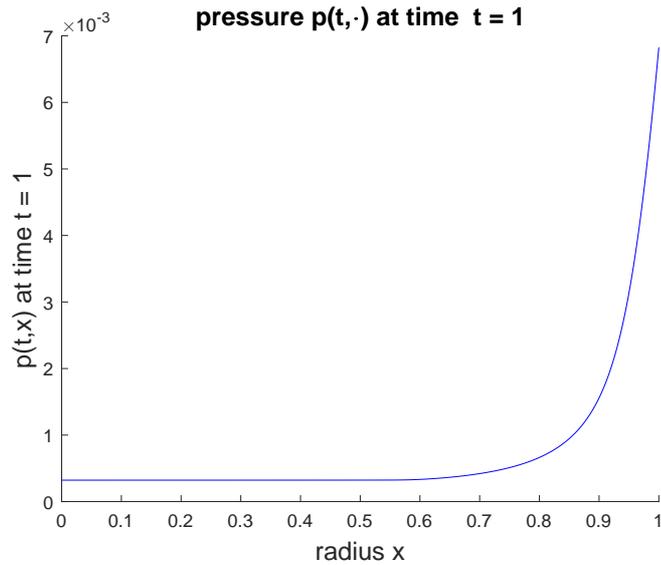}
\caption{Pressure $p$ at $t=1$ from the second example}
\label{ex2p}
\end{figure}
\end{center}
The numerical approximation leads us to the assumption that the exact solution depends only on $x/t$. Indeed, the existence of such a self-similar solution is justified in Lai's recent paper \cite[Theorem 1.1]{GL}.
Then we have a region emanating from the zero point
with a low constant pressure $p=0.00032$ and zero velocity $v=0$ for $t>0$,
followed by a centered rare\-fac\-tion fan starting from the zero point above the 
region with the constant values.
The numerical solution with $x_*=1$, $t_*=1$ and $N=3000$ is given in the Figures
\ref{ex2p} and \ref{ex2v}.
We also found that the computational values are in good agreement with those predicted from the results in Lai \cite{GL}.
\begin{center}
\begin{figure}[ht]
\includegraphics[width=0.75\textwidth]{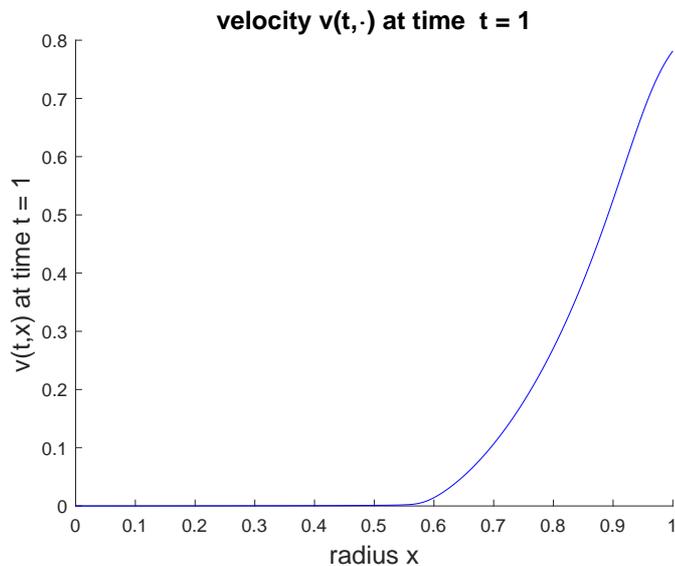}
\caption{Velocity $v$ at $t=1$ from the second example}
\label{ex2v}
\end{figure}
\end{center}
\item[3)] We choose the constant initial data $a_0=7$, $b_0=-4\sqrt{2}$ corresponding to a constant initial pressure $p_0=1$ and a constant radial part
$u_0=-1$ and $v_0=-1/\sqrt{2}$ of the initial four velo\-city and usual velo\-city,
respectively. The exact solution again depends
only on $x/t$, see \cite[Theorem 1.1]{GL}. 
Here we observe a straight line shock wave with slope $s=0.523$ emanating from the zero point, with a constant pressure $p=25.55$ and zero velocity $v=0$ for $t>0$ below the shock wave,
followed by a centered rare\-fac\-tion fan starting from the zero point above the shock wave.
The numerical solution with $x_*=1$, $t_*=1$ and $N=3000$ is given in Figures
\ref{ex3p} and \ref{ex3v}.
\begin{center}
\begin{figure}[ht]
\includegraphics[width=0.75\textwidth]{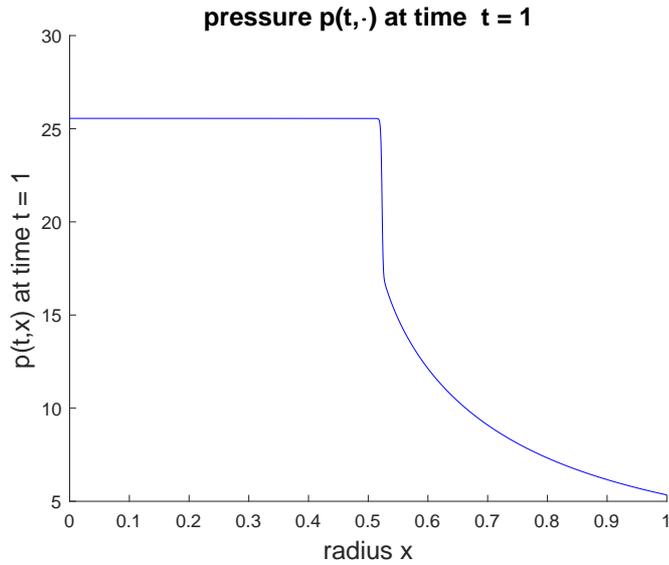}
\caption{Pressure $p$ at $t=1$ from the third example}
\label{ex3p}
\end{figure}
\end{center}
\begin{center}
\begin{figure}[ht]
\includegraphics[width=0.75\textwidth]{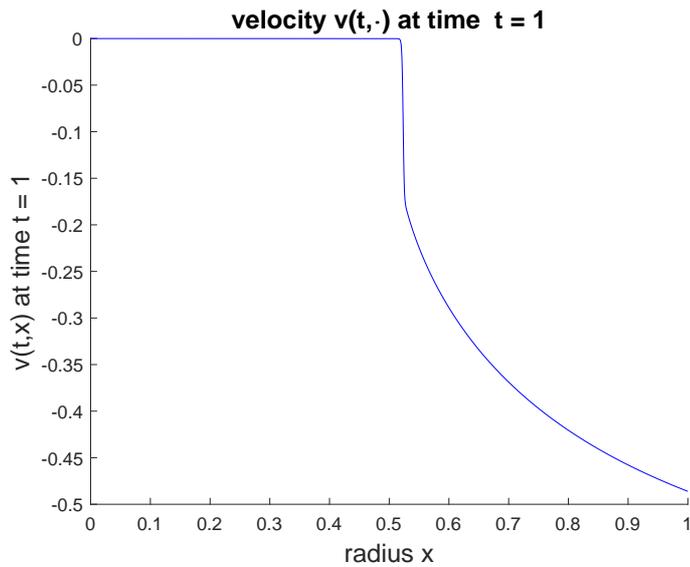}
\caption{Velocity $v$ at $t=1$ from the third example}
\label{ex3v}
\end{figure}
\end{center}
\item[4)] Expansion of a three dimensional spherical bubble with initial data
\begin{equation*}
p_0(x)=\begin{cases}1 \quad & \mbox{for~} 0 \leq x \leq 1\\
0.1 \quad & \mbox{for~} x > 1\,,\\
\end{cases}
\qquad
v_0(x)=0\,.
\end{equation*}

\begin{center}
\begin{figure}[ht]
\includegraphics[width=0.9\textwidth]{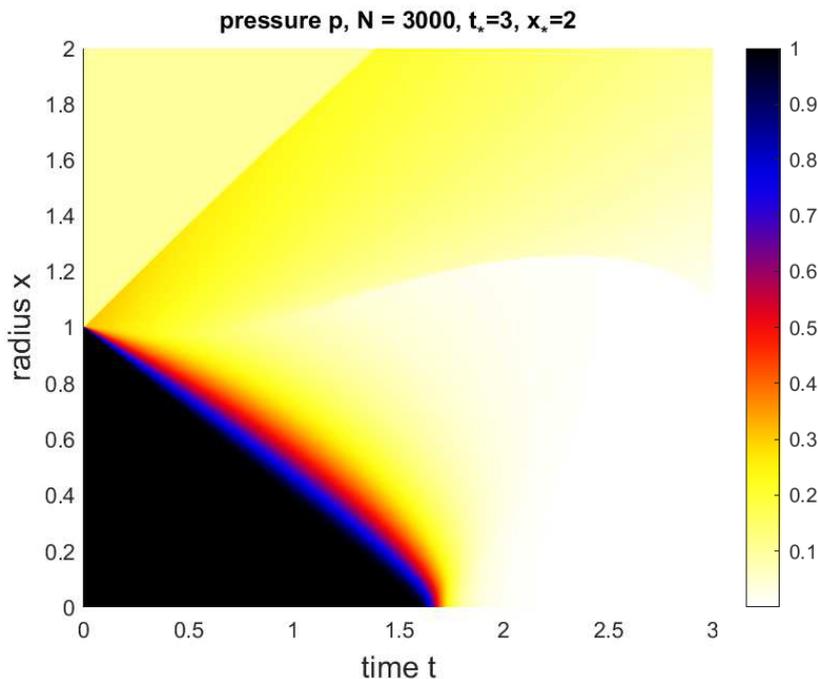}
\caption{Pressure $p$ from Example 4\,.}
\label{ex4p_color_time3}
\end{figure}
\end{center}

\begin{center}
\begin{figure}[ht]
\includegraphics[width=0.9\textwidth]{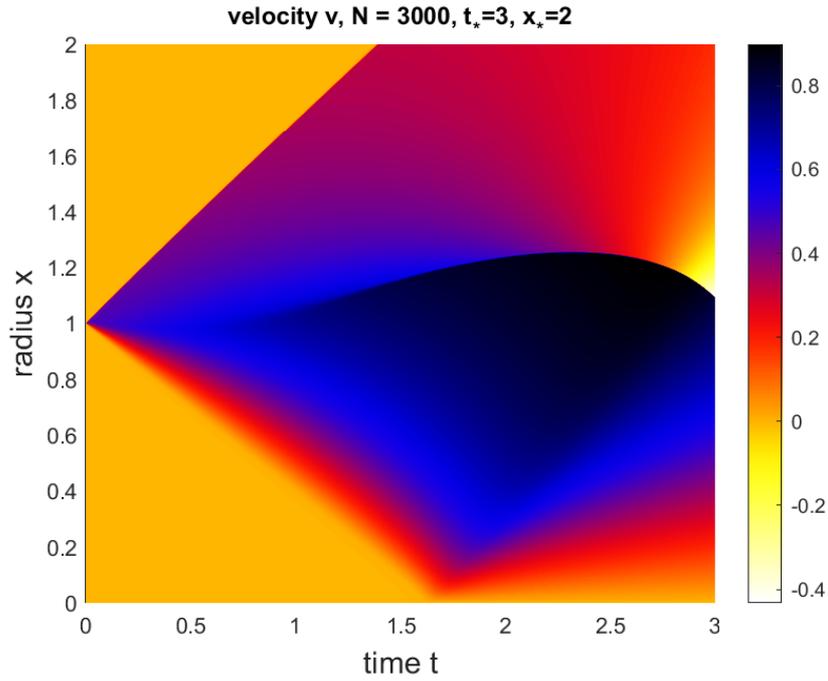}
\caption{Velocity $v$ from Example 4\,.}
\label{ex4v_color_time3}
\end{figure}
\end{center}
Initially, the pressure inside the bubble is ten times larger than outside,
which leads to a fast expansion of the bubble into the outer low pressure area.
This in turn gives rise to the formation of another 
low pressure area, 
namely the light yellow or white region in Figure \ref{ex4p_color_time3}
emanating from the zero point. 
The corresponding velocity is depicted in Figure \ref{ex4v_color_time3}.
We observe the formation of a shock wave,
running downwards into the new low pressure area
and reaching the zero point around time $t=4.16$, 
see Figures \ref{ex4p_color_time5} and \ref{ex4v_color_time5}.
The formation of this new shock wave is a peculiar nonlinear phenomenon. 
Shortly before the shock reaches the zero point the pressure takes very low values,
but its reflection from the zero point causes
a strong increase of the pressure in a very small time-space range near the boundary. 

\begin{center}
\begin{figure}[ht]
\includegraphics[width=0.9\textwidth]{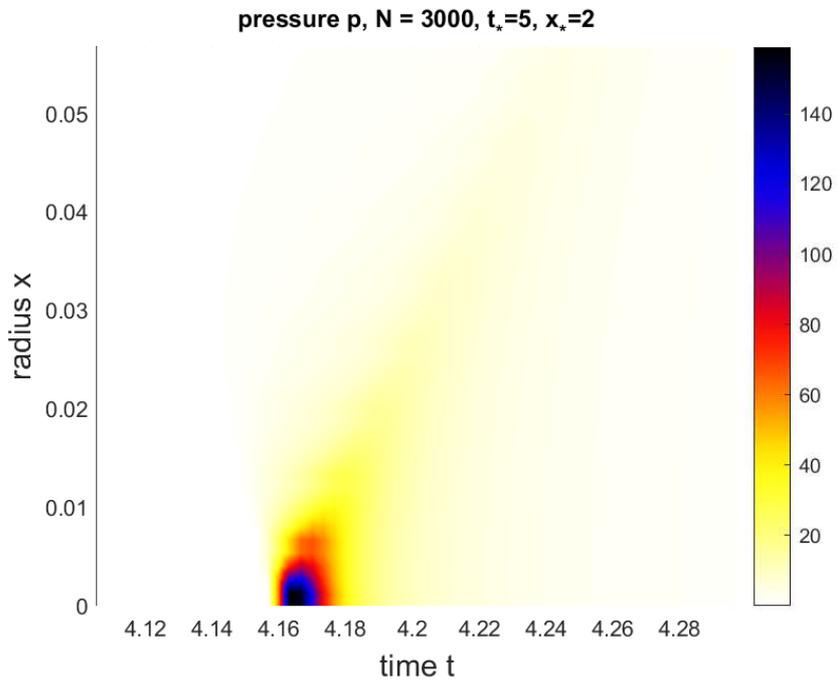}
\caption{Zoom near the pressure singularity, Example 4\,.}
\label{ex4p_color_time5}
\end{figure}
\end{center}

\begin{center}
\begin{figure}[ht]
\includegraphics[width=0.9\textwidth]{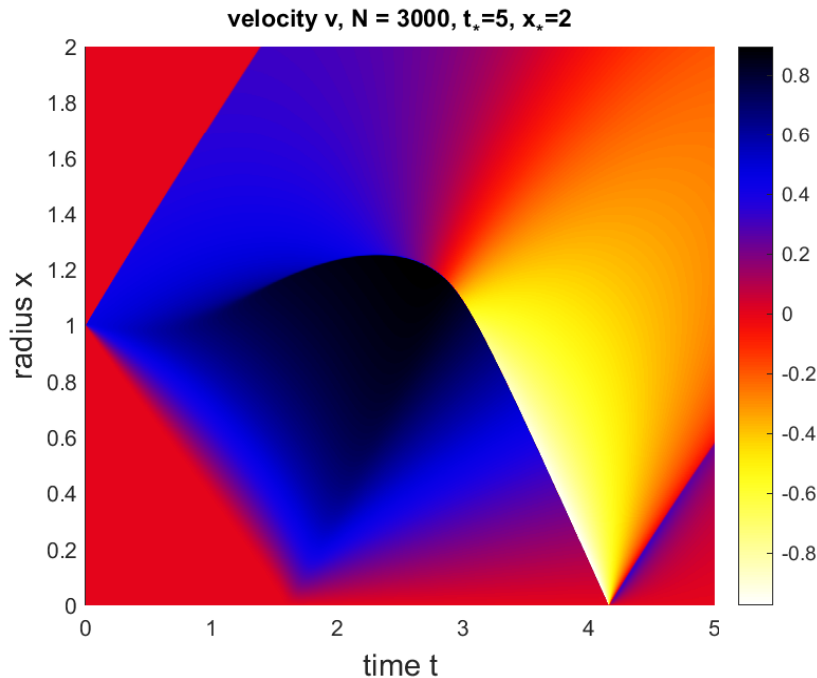}
\caption{Velocity $v$ from Example 4\,,
extension of the solution from Figure \ref{ex4v_color_time3}
with rescaled colors\,.}
\label{ex4v_color_time5}
\end{figure}
\end{center}

For the last example we have also changed the size of the initial bubble.
We only obtained the expected numerical solutions which are rescaled versions of the solutions presented here. Hence it is sufficient to study the problem 
with the initial bubble in the unit sphere around the origin.\\

Example 4 shows considerable differences
in the values of the pressure, especially 
in the domain $(t,x) \in [4\,,\,4.1]\times [0\,,\,0.02]$
a pressure less than $10^{-5}$.
At present we restrict our study to weak solutions without a vacuum state $a=b=0$. 
But a vacuum state may occur for certain initial data 
with symmetry and positive pressure in radially symmetric solutions of
the ultra-relativistic Euler equations, see
Lai's paper \cite[Lemmas 2.4, 2.5; Remark 2.2]{GL}.
In this case it is convenient to use the original
quantities $a$ and $b$ for which we have developed 
the scheme in Section \ref{scheme1}.
The question arises whether our scheme has the capability 
to capture more general solutions including the vacuum state accurately.\\
\end{itemize}

%%%%%%%%%%%%%%%%%%%%%%%%
{\bf Acknowledgement:} We thank Christoph Matern and Alexander Kaina for their support
in solving problems with Latex.

\bibliographystyle{plain}

\end{document}